\documentclass[11pt,a4paper,twoside]{article}
\usepackage{enumerate}
\usepackage{amsmath}
\usepackage{fix-cm}
\usepackage{pb-diagram}
\usepackage[english]{babel}
\usepackage{amssymb,latexsym}
\usepackage{graphicx}
\usepackage{hyperref}
\usepackage{color}
\usepackage{authblk}
\usepackage{tikz}
\DeclareMathOperator{\ric}{Ric}

\DeclareMathOperator{\scal}{Scal}

\def \bui#1#2{\mathrel{\mathop{\kern 0pt#1}\limits^{#2}}}
\def \buil#1#2{\mathrel{\mathop{\kern 0pt#1}\limits_{#2}}}

\newcommand{\II}{\mathrm{II}}
\let\<\langle
\let\>\rangle
\newcommand{\R}{{\mathbb R}}

\newtheorem{example}{Examples}[section]
\newtheorem{thm}{Theorem}[section]
\newtheorem{lemma}[thm]{Lemma}
\newtheorem{prop}[thm]{Proposition}

\newtheorem{cor}[thm]{Corollary}
\newtheorem{remark}[thm]{Remark}
\newtheorem{remarks}[thm]{Remarks}
\newtheorem{definition}[thm]{Definition}
\newtheorem{notation}[thm]{Notation}
\newtheorem{exabout:ample}[thm]{Example}

\title{Integral Formulas for Differential Forms on Weighted Manifolds and Applications }

\author[1]{Fida El Chami\thanks{\texttt{fchami@ul.edu.lb}}}
\author[1]{Ola Makhoul\thanks{\texttt{ola.makhoul@ul.edu.lb}}}

\affil[1]{\footnotesize Lebanese University, Faculty of Sciences II, Department of Mathematics, P.O. Box 90656 Fanar-Matn, Lebanon}

\begin{document}
\maketitle


\noindent\begin{center}\begin{tabular}{p{115mm}}
\begin{small}{\bf Abstract.}
In this paper, we derive a Reilly formula for differential forms on weighted manifolds with nonempty boundary. As an application of this formula, we prove a Poincar\'e-type inequality in the same context and explore several of its consequences. We also present weighted versions of some boundary value problems and obtain new eigenvalue estimates that extend previously known results.

\end{small}\\
\end{tabular}\end{center}

\noindent\begin{small}{\it Mathematics Subject Classification} (2020): 53C21, 58C40, 58J32, 58J50
\end{small}

\noindent\begin{small}{\it Keywords}: weighted manifolds with boundary, Reilly formula, Poincar\'e inequality, boundary value problems, eigenvalue estimates  	
\end{small}

\section{Introduction} 
The study of geometric inequalities on manifolds with boundary often begins with powerful integral formulas. One important example is the Reilly formula, established by Robert Reilly \cite{Reilly1977}, which applies to functions defined on an $n$-dimensional compact Riemannian manifold $(M,g)$ with smooth boundary $\partial M$. This formula provides a powerful tool for geometric analysis on manifolds with boundary and serves as a fundamental tool that connects curvature, boundary geometry, and analytic properties of functions, and has been extensively used to derive many important estimates, particularly for the eigenvalues of various operators as well as rigidity results. Later, Miao and Wang \cite{MiaoWang2016} used this formula to prove a sharp Poincaré-type inequality for functions under a lower bound on the Ricci curvature. This inequality has led to important geometric comparison results and rigidity theorems including that equality holds precisely when the manifold is isometric to a Euclidean ball. \\\\
In 2011, Raulot and Savo \cite{RaulotSavo2011} established a generalization of Reilly's formula to differential forms. Their formula states that for all $\omega \in \Omega^p(M)$:
\begin{eqnarray}\label{eq:Reillyformula}
\int_M\left(|d\omega|^2+|\delta\omega|^2\right)\, d\mu_g=\int_M\left(|\nabla\omega|^2+\langle W^{[p]}\omega,\omega\rangle\right) \,d\mu_g+\int_{\partial M}\left( 2\langle\nu\lrcorner\omega,\delta^{\partial M}(\iota^*\omega)\rangle+\mathcal{B}(\omega,\omega)\right)\, d\mu_g, 
 \end{eqnarray}
 where $W^{[p]}$ is the curvature operator in the Weitzenb\"ock formula and
 \begin{eqnarray*}
 	\mathcal{B}(\omega,\omega)=\langle S^{[p]}(\iota^*\omega),\iota^*\omega\rangle+(n-1)H|\nu\lrcorner\omega|^2-\langle S^{[p-1]}(\nu\lrcorner\omega),\nu\lrcorner\omega\rangle, 
 \end{eqnarray*}
  $\nu$ being the inward unit vector field normal to $\partial M$, $\delta^{\partial M}$ the codifferential on $\partial M$ and $S^{[p]}$ the extension of the shape operator to differential $p$-forms.
  This formula has proven to have several important applications. Raulot and Savo, for instance, used it to obtain a sharp lower bound for the spectrum of the boundary Hodge Laplacian and, in the same work, established several rigidity results. Subsequently, Ginoux, Habib and Raulot \cite{GHR} used this formula to establish an extension of the so-called Poincar\'e inequality to differential forms. Specifically, they proved that for any exact $p$-form $\omega$ on $\partial M$, 
 \begin{equation} \label{Poincare1}
 \displaystyle \int_{\partial M} \langle S^{[p]} \omega, \omega \rangle d\mu_g \le \int_{\partial M} \frac{|\delta^{\partial M}w|^2}{\sigma_{n-p}} d\mu_g, \end{equation}
 provided $W^{[p]} \ge 0$ and $\sigma_{n-p} >0$, where $\sigma_{n-p}$ denotes the sum of the  $(n-p)$ smallest principal curvatures of $\partial M$. Another version of this inequality was also obtained under the assumption that the curvature term $W^{[p]}$ is bounded below by a positive constant.
 As an application, they derived a new inequality involving the mean and scalar curvatures of the boundary and characterized its limiting case in codimension one and thus extending the rigidity results of Miao and Wang \cite{MiaoWang2016}. \\
In another direction, in $\cite{EGHMR}$, El Chami, Ginoux, Habib, Makhoul and Raulot introduced the buckling problem on differential forms
 \begin{equation}\label{buckling_forms}
 \left\{
 \begin{array}{lll}
 \Delta^2 \omega&=\Lambda \Delta \omega &\textrm{ on } M\\ \omega
 &=0&\textrm{ on }\partial M\\
 \nabla_\nu\omega&=0&\textrm{ on }\partial M \end{array}\right.
 \end{equation}
 for some real constant $\Lambda$. In the same paper, they defined the clamped plate problem also on differential forms
 \begin{equation}\label{clamped_plate_forms}
 \left\{
 \begin{array}{lll}\Delta^2 \omega&=\Gamma \omega &\textrm{ on } M\\ 
 \omega &=0&\textrm{ on }\partial M\\
 \nabla_\nu\omega&=0&\textrm{ on }\partial M 
 \end{array}
 \right.
 \end{equation}
 for some real constant $\Gamma$. 
 For both problems, which generalize classical fourth-order problems from functions to forms, they studied their spectral properties, characterized their first eigenvalues, and obtained as well some estimates involving it, using the Reilly formula (\ref{eq:Reillyformula}).
\\\\
In a different setting, in 2010, Ma and Du \cite{MaDu2010} also generalize the (scalar) Reilly formula to weighted compact manifolds with boundary, obtaining bounds for the first nonzero eigenvalue of the weighted Laplacian under both Dirichlet and Neumann boundary conditions.   
Afterwards, a number of works have further developed this line of research on weighted manifolds. For instance, Huang and Ruan \cite{HuangRuan2014} established estimates for the weighted mean curvature and used them to derive new lower bounds for eigenvalues of the weighted Laplacian. Later, Batista, Cavalcante, and Pyo \cite{BCP2014} studied compact submanifolds of weighted manifolds, obtaining isoperimetric inequalities and related eigenvalue estimates under curvature assumptions. Tu and Hang \cite{TuHuang2019} used the weighted Reilly formula to study the effect of the $m$-dimensional Bakry--\'Emery Ricci curvature on the geometry of the boundary of $M$.
In 2019, Ilias and Shouman \cite{IliasShouman2019} gave comparison results for the eigenvalues of various eigenvalue problems, including Dirichlet, Neumann, clamped, and buckling, on weighted manifolds under curvature assumptions. See also the work of Bezerra and Xia \cite{BezerraXia2020} on similar operators.
\\\\
In the present paper, we generalize the aforementioned results to the case of differential forms on weighted manifolds with boundary. 
Before outlining the contents of this paper, let us mention some additional  references directly relevant to our discussion on differential forms on weighted manifolds. 
In 2020, Petersen and Wink \cite{PetersenWink20} studied the Bochner formula for the weighted Laplacian acting on differential forms. Later, in 2024, Branding and Habib \cite{BrandingHabib22} obtained some estimates concerning the eigenvalues of the Hodge Laplacian acting on differential forms on a Bakry--\'Emery manifold without boundary.  
\\\\ 
The paper is organized as follows. In section $2$, we start with some useful preliminaries concerning compact manifolds with smooth boundary, endowed with a weighted measure. Then in section $3$, we generalize the Reilly formula \eqref{eq:Reillyformula} to the context of weighted manifolds (see Theorem \ref{weightedReillyformula}).  In section $4$, we briefly present the extension of the cohomology of compact Riemannian manifolds with boundary to the same context, based on the study of Schwarz \cite{Sc}. Note that the weighted relative and absolute cohomologies turn out to be isomorphic to the unweighted ones. The generalized Reilly formula \eqref{WReilly_formula}, with the help of the abovementioned isomorphism, allows us to establish a weighted version of the Poincar\'e inequality \eqref{Poincare1} in the case where the $p$-Ricci tensor is nonnegative and in the case where it is bounded below by a real and nonnegative number. These Poincar\'e inequalities are then used to derive new results. For instance, we obtain upper bounds for the integral of the mean curvature with respect to the weighted measure in terms of quantities involving the weight function $f$ and some curvature terms of $\partial M$ in the case where the boundary of the manifold is isometrically immersed in Euclidean space, and when it is immersed in the sphere. Finally in section $5$, we give a generalization of the buckling problem \eqref{buckling_forms} 
and the clamped plate problem \eqref{clamped_plate_forms} 
to the same framework. Moreover, we extend some results related to the smallest eigenvalues of the mentioned eigenvalue problems. \\\\
\textbf{Acknowledgement}: We would like to thank N.~Ginoux, G.~Habib and S.~Raulot for their helpful comments and constructive suggestions regarding this work. This work was supported by a grant from the Lebanese University.
\section{Preliminaries}
Let $(M^n,g)$ be a compact Riemannian manifold with smooth boundary. A weighted Riemannian manifold $(M,g, d\mu_f)$ (also called a Bakry--\'Emery manifold) is a Riemannian manifold $(M,g)$ endowed with a weighted measure $d \mu_f= e^{-f} d\mu_g$, where $f$ is a real-valued smooth function on $M$ and $d \mu_g$ is the Riemannian measure
induced by the metric $g$. The Bakry--\'Emery Laplacian, also called the weighted Laplacian, on differential forms on $M$ is defined to be $\Delta_f:= d \delta_f +\delta_f d$, where $d$ is the exterior derivative
and $\delta_f=\delta + \nabla f \lrcorner \, $ (see for instance \cite{PetersenWink20} and \cite{BrandingHabib22}). It is easy to see that the weighted Laplacian $\Delta_f$ can be computed using the Hodge Laplacian by the formula $\Delta_f=\Delta + \mathcal{L}_{\nabla f}$, where $\mathcal{L}_{\nabla f}$ is the Lie derivative with respect to the gradient of $f$. Note that for any differential form $\omega$, $\delta_f \omega=e^{f} \delta (e^{-f} \omega)$ from which it can be deduced that $\delta_f^2=0$.

Recall the following Bochner formula in weighted manifolds (see e.g. \cite{PetersenWink20})
\begin{equation}
\label{WBochner}
\Delta_f =\nabla^*_f \nabla + \ric_f^{(p)},
\end{equation}
where $\nabla_f^*$ denotes the formal adjoint of $\nabla$ with respect to the weighted measure $d \mu_f$, and $\ric_f^{(p)}:=W^{[p]} + T^{[p]}_f$ is the $p$-Ricci tensor. Here $W^{[p]}$ is the Bochner operator 
$W^{[p]}:=-\sum_{i,j=1}^n e_i^* \wedge(e_j\lrcorner \, R(e_i,e_j))$, 
where $\{e_i\}_{i=1, \dots ,n}$ is
a local orthonormal frame of $TM$ and $R$ is the Riemann curvature tensor of $(M,g)$ given by $R(X,Y)=[\nabla_X, \nabla_Y]-\nabla_{[X,Y]}$  for all $X,Y$ tangent to $M$,  
and $T^{[p]}_f$ is the self adjoint endomorphism of $\Omega^p(M)$ given by
$$T^{[p]}_f (\omega)(X_1, \dots , X_p):=\sum_{j=1}^p \omega(X_1, \dots , \nabla^2 f (X_j), \dots , X_p),$$
for any $\omega \in \Omega^p(M)$ and $X_1, \dots , X_p \in \Gamma(TM)$. Note that for $p=1$, $\ric_f^{(1)}$ is the $\infty$-Bakry--\'Emery tensor introduced by Lichnerowicz \cite{Lichnerowicz1, Lichnerowicz2} and by Bakry and \'Emery \cite{Bakry-Emery}.
Note that $\nabla_f^* \nabla=\nabla^*\nabla +\nabla_{\nabla f}$.

Next, we recall some additional notions that will be useful later. Let $S$ be the shape operator given by $S(X)=-\nabla_X \nu$ for $X \in \Gamma(T\partial M)$, $H=\frac{1}{n-1} \textrm{tr}(S)$ the mean curvature of the boundary and $S^{[p]}$ the canonical extension of $S$ to differential forms defined as follows:
$$S^{[p]}(\omega)(X_1, \dots , X_p)=\sum_{j=1}^p \omega (X_1, \dots , S(X_j), \dots , X_p),$$
for $\omega \in \Omega^p(\partial M)$ and $X_1, \dots, X_p \in \Gamma(T\partial M)$. For $p=0$, by convention, $S^{[0]}=0$. We denote by $\eta_1, \dots , \eta_{n-1}$ the principal curvatures of $\partial M$ arranged so that $\eta_1 \le \dots \le \eta_{n-1}$. For $p \in \{1, \dots , n-1\}$, let $\sigma_p=\eta_1 + \dots + \eta_p$ be the lowest $p$-curvature, and let $\sigma_p(\partial M)=\displaystyle\inf_{x \in \partial M} \sigma_p(x)$. One can easily check that $$\langle S^{[p]}(\omega), \omega \rangle \ge \sigma_p  |\omega|^2.$$

\section{A Reilly formula for differential forms on weighted manifolds}
In this section, we provide a weighted version of the Reilly formula for differential forms \eqref{eq:Reillyformula}. This result generalizes \cite[Thm. 3]{RaulotSavo2011} to the case where the Riemannian manifold $(M,g)$ is equipped with a weighted measure $d\mu_f$, $f$ being a smooth function on $M$.
As a first step toward proving the Reilly formula in weighted manifolds, we give a weighted version of the partial integration formula.
\begin{lemma}
    Let $(M^n,g,d\mu_f)$ be a compact weighted manifold with smooth boundary $\partial M$ and let $\nu$ be the inward unit normal vector field along the boundary.  For any $\alpha \in \Omega^p(M)$, $\beta \in \Omega^{p+1}(M)$, we have
    \begin{eqnarray}
\label{stokes_drift}
    \int_M \langle d\alpha, \beta \rangle \, d \mu_f=\int_M \langle \alpha, \delta_f \beta \rangle \, d \mu_f -\int_{\partial M} \langle \iota^* \alpha, \nu \lrcorner \, \, \beta \rangle \, d \mu_f.
\end{eqnarray}
\end{lemma}
{\it Proof.}
Let $\alpha \in \Omega^p(M)$ and $\beta \in \Omega^{p+1}(M)$. Applying the integration by parts formula to $\alpha$ and $e^{-f} \beta$ gives
$$\int_M \langle d\alpha, e^{-f}\beta \rangle \, d \mu_g \nonumber\\
     =\int_M \langle \alpha, \delta (e^{-f}\beta) \rangle \, d \mu_g -\int_{\partial M} \langle \iota^* \alpha, \nu \lrcorner \, \, (e^{-f}\beta) \rangle \, d \mu_g ,$$
     so that
\begin{eqnarray*}
    \int_M \langle d\alpha, \beta \rangle \, d \mu_f 
     &=& \int_M \langle \alpha, \delta \beta \rangle \, d \mu_f + \int_M \langle \alpha, \nabla f \lrcorner \, \, \beta \rangle \, d\mu_f -\int_{\partial M} \langle \iota^* \alpha, \nu \lrcorner \, \, \beta \rangle \, d \mu_f \nonumber \\
      &=& \int_M \langle \alpha, \delta_f \beta \rangle \, d \mu_f -\int_{\partial M} \langle \iota^* \alpha, \nu \lrcorner \, \beta \rangle \, d \mu_f . 
\end{eqnarray*} \hfill$\square$ \\ \\
As a consequence, we obtain the following formula:
\begin{lemma}
    Let $\omega$ be a $p$-form on $M$. Then
\begin{equation}
    \label{Deltaomom_drift}
\int_M |d \omega|^2 \, d \mu_f +\int_M |\delta_f \omega|^2 \, d \mu_f = \int_M\langle \Delta_f \omega , \omega \rangle  \, d \mu_f+ \int_{\partial M} \langle \iota^* \delta_f \omega , \nu \lrcorner \, \omega \rangle \, d \mu_f -\int_{\partial M} \langle \nu \lrcorner \, d \omega , \iota^* \omega \rangle \, d \mu_f. 
\end{equation}
\end{lemma}
{\it Proof.}
For any forms $\omega, \omega' \in \Omega^p(M)$, a straightforward computation, using \eqref{stokes_drift}, yields:
  \begin{equation*}
\int_M\langle\Delta_f\omega,\omega'\rangle \, d \mu_f
=\int_M\langle d\omega,d\omega'\rangle \, d \mu_f + \int_M\langle\delta_f\omega,\delta_f\omega'\rangle \, d \mu_f-\int_{\partial
M}\langle \iota^*\delta_f\omega,\nu\lrcorner \,\,\omega'\rangle \, d \mu_f
+\int_{\partial
M}\langle\nu\lrcorner \, d\omega,\iota^*\omega'\rangle \, d \mu_f.
\end{equation*}
The statement follows by taking $\omega'=\omega$. \hfill$\square$ 
\\ \\
Let us recall the following formulas \cite[Lem. 18]{RaulotSavo2011} valid in a Riemannian manifold with boundary: for any differential form $\omega$ on $M$, we have: 
 \begin{eqnarray}
 \label{ddeltabord}
\left\{
  \begin{array}{lll}
  \medskip
  \delta^{\partial M} (\iota^* \omega)&=& \iota^*(\delta \omega)+\nu \lrcorner \, \nabla_\nu \omega +S^{[p-1]}(\nu \lrcorner \, \omega) -(n-1)H \nu \lrcorner \, \omega\\
  \medskip
  d^{\partial M} (\nu \lrcorner \, \omega)&=& -\nu \lrcorner \, d \omega+\iota^*( \nabla_\nu \omega) -S^{[p]}(\iota^* \omega),
  \end{array}
\right.
 \end{eqnarray}
where $\iota : \partial M \longrightarrow M$ is the canonical embedding, and $d^{\partial M}$ and $ \delta^{\partial M}$ stand, respectively, for the exterior derivative and codifferential on the boundary. \\ \\
Since $\iota^*(\nabla f \lrcorner\, \omega) = \nabla^{\partial M} f  \lrcorner \,\iota^*\omega + f_\nu \, \nu \lrcorner \, \omega$ for $\omega \in \Omega^p(M)$, where $f_\nu:=\partial_\nu f$, we can directly extend the first identity of \eqref{ddeltabord} to the case of weighted manifolds with boundary. Namely, we obtain:
\begin{lemma}
      Let $(M^n,g,d\mu_f)$ be a compact weighted manifold with smooth boundary $\partial M$. For any differential form $\omega \in \Omega^p(M)$, we obtain:
      \begin{equation}
 \label{deltabord_Weighted}
  \delta_f^{\partial M} (\iota^* \omega)= \iota^*(\delta_f \omega)+\nu \lrcorner \, \nabla_\nu \omega +S^{[p-1]}(\nu \lrcorner \, \omega) -(n-1)H_f \, \nu \lrcorner \, \omega.
 \end{equation}
Here, $\delta_f^{\partial M}$ refers to the codifferential on the weighted manifold $\partial M$ and $H_f:=H+\dfrac{1}{n-1}f_\nu$ is the weighted mean curvature. 
\end{lemma}
We are now able to prove the main result of this section.
 \begin{thm} [Weighted Reilly formula]\label{weightedReillyformula} 
     Let $(M^n,g,d\mu_f)$ be a compact weighted manifold with smooth boundary $\partial M$. For any $p$-form $\omega$ on $M$, we have 
     \begin{eqnarray}
         \label{WReilly_formula}
          \int_M |d \omega|^2 \, d \mu_f +\int_M |\delta_f \omega|^2 \, d \mu_f &=& \int_M |\nabla \omega|^2 \, d \mu_f +\int_M \langle \ric_f^{(p)} \omega , \omega \rangle \, d\mu_f
   \nonumber \\&&
     + 2\int_{\partial M} \langle \delta^{\partial M}_f (\iota^*\omega) , \nu \lrcorner \, \omega \rangle \, d \mu_f  + \int_{\partial M} \mathcal{B}_f(\omega, \omega) \, d \mu_f,
       \end{eqnarray}
where $\mathcal{B}_f(\omega, \omega)= \langle S^{[p]}(\iota^* \omega), \iota^* \omega \rangle-\langle S^{[p-1]}(\nu \lrcorner \, \omega), \nu \lrcorner \, \omega \rangle +(n-1)  H_f \, |  \nu \lrcorner \, \omega |^2$.
 \end{thm}
 {\it Proof.} 
We proceed as in the proof of \cite[Thm. 3]{RaulotSavo2011}. Let $\omega$ be a $p$-form on $M$. 
First, since $ \nabla^*_f \nabla =\nabla^*\nabla  +\nabla_{\nabla f}$, we obtain
\begin{eqnarray*}
    \langle \nabla^*_f \nabla \omega, \omega \rangle &=& |\nabla \omega|^2 + \dfrac{1}{2}\Delta(|\omega|^2) + \langle \nabla_{\nabla f} \omega, \omega \rangle \\
     &=& |\nabla \omega|^2 + \dfrac{1}{2}\Delta(|\omega|^2) +\dfrac{1}{2} \langle d f ,d|\omega|^2\rangle \\
     &=& |\nabla \omega|^2 + \dfrac{1}{2}\Delta_f(|\omega|^2).
\end{eqnarray*}
Then, by integrating this last identity over $M$ with respect to the measure $d \mu_f$, we get
\begin{eqnarray*}
    \int_M \langle \nabla^*_f \nabla \omega, \omega \rangle \, d\mu_f&=&\int_M |\nabla \omega|^2 \, d \mu_f +\dfrac{1}{2}\int_M \Delta_f(| \omega|^2) \, d \mu_f \\
    &=& \int_M |\nabla \omega|^2 \, d \mu_f +\int_{\partial M} \langle \nabla_\nu \omega, \omega\rangle \, d \mu_f \\{\Huge }
    &=& \int_M |\nabla \omega|^2 \, d \mu_f +\int_{\partial M} \langle \iota^* \nabla_\nu \omega, \iota^*\omega\rangle \, d \mu_f + \int_{\partial M} \langle \nu \lrcorner \, \nabla_\nu \omega, \nu \lrcorner \,\omega\rangle \, d \mu_f .
\end{eqnarray*}
On the other hand, we integrate the Bochner formula \eqref{WBochner} and substitute the term $\int_M \langle \nabla^*_f \nabla \omega, \omega \rangle \, d\mu_f$ from the previous equality, yielding
\begin{eqnarray*}
    \int_M \langle \Delta_f \omega , \omega \rangle \, d\mu_f 
    &=& \int_M \langle \nabla^*_f \nabla \omega, \omega \rangle \, d\mu_f +\int_M \langle \ric_f^{(p)} \omega , \omega \rangle \, d\mu_f \nonumber \\
    &=& \int_M |\nabla \omega|^2 \, d \mu_f+  \int_M \langle \ric_f^{(p)} \omega , \omega \rangle\, d \mu_f \nonumber \\
    && +\int_{\partial M} \langle \iota^* \nabla_\nu \omega, \iota^*\omega\rangle \, d \mu_f + \langle \nu \lrcorner \, \nabla_\nu \omega, \nu \lrcorner \,\omega\rangle\, d \mu_f.
\end{eqnarray*}

Thus, combining this equality with \eqref{Deltaomom_drift}, we deduce
\begin{eqnarray*}
    \int_M |d \omega|^2 \, d \mu_f +\int_M |\delta_f \omega|^2 \, d \mu_f &=& \int_M |\nabla \omega|^2 \, d \mu_f +\int_M \langle \ric_f^{(p)} \omega , \omega \rangle \, d\mu_f\\
    &+&\int_{\partial M} \langle \iota^* \nabla_\nu \omega, \iota^*\omega\rangle \, d \mu_f + \int_{\partial M} \langle \nu \lrcorner \, \nabla_\nu \omega, \nu \lrcorner \,\omega\rangle \, d \mu_f 
    \\
   &&  \int_{\partial M} \langle \iota^* \delta_f \omega , \nu \lrcorner \, \omega \rangle \, d \mu_f -\int_{\partial M} \langle \nu \lrcorner \, d \omega , \iota^* \omega \rangle \, d \mu_f \\
    &=& \int_M |\nabla \omega|^2 \, d \mu_f  +\int_M \langle \ric_f^{(p)} \omega , \omega \rangle \, d\mu_f +\int_{\partial M} \langle  \delta^{\partial M}_f (\iota^*\omega) , \nu \lrcorner \, \omega \rangle \, d \mu_f \\ && - \int_{\partial M} \langle
 S^{[p-1]}(\nu \lrcorner \, \omega), \nu \lrcorner \, \omega \rangle \, d \mu_f + (n-1)\int_{\partial M} H_f \, |  \nu \lrcorner \, \omega |^2 \, d \mu_f \\ 
 && + \int_{\partial M} \langle d^{\partial M}(\nu \lrcorner \, \omega) , \iota^* \omega \rangle \, d \mu_f +\int_{\partial M} \langle S^{[p]}(\iota^* \omega), \iota^* \omega \rangle \, d \mu_f.
\end{eqnarray*}
In the last equality, we used Equations  \eqref{ddeltabord} and \eqref{deltabord_Weighted}. By observing that $\int_{\partial M} \langle d^{\partial M}(\nu \lrcorner \, \omega) , \iota^* \omega \rangle \, d \mu_f =\int_{\partial M} \langle \nu \lrcorner \, \omega , \delta_f^{\partial M}(\iota^* \omega) \rangle \, d \mu_f $, we obtain the announced formula \eqref{WReilly_formula}.
\hfill$\square$ 
\begin{remark}
    Based on the identity $\star_{\partial M}S^{[p]}+S^{[n-1-p]}\star_{\partial M}=(n-1)H \star _{\partial M}$ on $p$-forms of $\partial M$, and due to the fact that $\iota^* \star \omega$ is equal to $\star_{\partial M}(\nu \lrcorner \, \omega)$ up to sign, the boundary term $\mathcal{B}_f(\omega, \omega)$ can also be expressed the following way:
    $$\mathcal{B}_f(\omega, \omega)= \langle S^{[p]}(\iota^* \omega), \iota^* \omega \rangle +\langle S^{[n-p]}(\iota^* \star \omega), \iota^* \star \omega \rangle +f_\nu | \nu \lrcorner \omega |^2.$$ 
\end{remark}
\begin{remark} 
\begin{enumerate}
    \item In the special case where $p=1$, $\omega=du$, $u$ being a smooth function, the Reilly formula \eqref{WReilly_formula} reduces to 
    \begin{eqnarray*}
          \int_M (\Delta_f u)^2 \, d \mu_f &=& \int_M \big(|\nabla^2 u|^2 + \ric_f (\nabla u , \nabla u) \big) \, d\mu_f
   \nonumber \\&&
     + \int_{\partial M}  \big( 2\, (\Delta^{\partial M}_f u ) u_\nu \, + \langle S(\nabla^{\partial M} u), \nabla^{\partial M} u \rangle + (n-1) H_f \, u_\nu \big)\, d \mu_f,
       \end{eqnarray*}
which is the weighted Reilly formula for functions established by Ma and Du in \cite[Thm. 1]{MaDu2010}.
    \item 
    In \cite[Thm. 1]{Xiong23}, Xiong also provided a weighted version of the Reilly formula for differential forms that generalizes \cite[Thm. 3]{RaulotSavo2011}. However, his result differs from ours, given that the context is not the same: we are working in a Bakry--\'Emery manifold with the weighted Laplacian, whereas Xiong's work is in the Riemannian setting and is based on the integration of a weighted version of the Bochner formula.
    
\end{enumerate}
\end{remark}

\section{A weighted Poincar\'e Inequality for differential forms and applications}
In this section, we establish a Poincar\'e-type inequality for differential forms in Bakry--\'Emery manifolds, under two different curvature assumptions and we then provide several applications. Our results generalize those obtained in \cite[section 3]{GHR} for differential forms on Riemannian manifolds. Note that in the scalar case, the same type of result was obtained by Tu and Huang \cite{TuHuang2019} in Bakry--\'Emery manifolds under the assumption that the $m$- dimensional Bakry--\'Emery Ricci curvature is bounded from below. 
\\\\
Before presenting our first main result, we will explain how to extend some facts concerning the cohomology of compact Riemannian manifolds with boundary to the context of weighted manifolds. These results are essentially based on the developments presented in Schwarz \cite{Sc} and can be extended to our framework without significant difficulty. \\
Let $(M,g,d\mu_f)$ be a compact weighted manifold with smooth boundary $\partial M$ and let $\nu$ be the inward unit normal vector field along $\partial M$. 
We denote by ${H}^1 \Omega^p_{D}(M):=\{\omega \in H^1 \Omega^p(M)\,|\, \iota^*\omega=0\}$ (resp. ${H}^1 \Omega^p_{N}(M):=\{\omega \in H^1\Omega^p (M) \,| \,  \nu \lrcorner\omega=0\}$) the spaces of $p$-forms of Sobolev class $H^1$ with respect to the weighted measure $d\mu_f$, with vanishing tangential (resp. normal components).
Also, let us denote by $\mathcal{H}^p_{f}(M):=\{\omega \in H^1 \Omega^p(M)\,|\, d\omega=\delta_f \omega =0\}$ the space of $f$-harmonic fields. 
The spaces $\mathcal{H}^p_{f,D}(M):=\{\omega \in H^1 \Omega^p(M)\,|\, d\omega=\delta_f \omega =0, \iota^*\omega=0\}$ and $\mathcal{H}^p_{f,N}(M):=\{\omega \in H^1\Omega^p (M) \,| \, d\omega=\delta_f \omega =0, \nu \lrcorner\omega=0\}$ are called the spaces of Dirichlet and Neumann $f$-fields, respectively. 
One proves exactly as in the classical case that (see \cite[Sections 2.2 \& 3.4]{Sc}) the spaces $\mathcal{H}^p_{f,D}(M)$ and $\mathcal{H}^p_{f,N}(M)$ are finite-dimensional consisting of smooth forms and that: 
      \begin{equation}
      \label{f-harmonic_fields_DN}
    \mathcal{H}^p_{f,D}(M) \cap \mathcal{H}^p_{f,N}(M)=\{0\}.
    \end{equation}
The Hodge-Morrey decomposition remains valid in a Bakry--\'Emery manifold: The Hilbert space $L^2 \Omega^p(M)$ splits into the $L^2$-orthogonal direct sum (with respect to the measure $d\mu_f$):
    $$L^2 \Omega^p(M) = \mathcal{E}^p(M) \oplus \mathcal{C}_f^p(M) \oplus L^2\mathcal{H}^p_f(M).$$
    Here, $\mathcal{E}^p(M):=\{d\alpha \, | \, \alpha\in H^1_D\Omega^{p-1}(M)\} $, 
$\mathcal{C}_f^p(M):=\{\delta_f\alpha \, | \, \alpha\in H^1_N\Omega^{p-1}(M)\}$, and $L^2\mathcal{H}^p_f(M):=\overline{\mathcal{H}^p_f(M)}$ is the $L^2$-closure of the space of $f$-harmonic fields. Furthermore, the Friedrichs decomposition also holds in this setting and its proof follows exactly the argument given in Schwarz \cite[Thm. 2.4.8]{Sc}. By using the Hodge-Morrey decomposition together with the Friedrichs decomposition, we can prove a variant of the Hodge isomorphism theorem for weighted manifolds with boundary, and establish an isomorphism between the absolute (resp. relative) weighted cohomology, the de Rham cohomology $H^p(M,d)$ (resp. $H^p(M,\delta_f)$)  and the harmonic cohomology with Neumann (resp. Dirichlet) boundary condition. In other words, the relative weighted cohomology and the absolute weighted cohomology are independent of $f$ in the sense that 
$$
\mathcal{H}^p_{f,N}(M) \cong \mathcal{H}^p_{N}(M) \quad \text{and} \quad \mathcal{H}^p_{f,D}(M) \cong \mathcal{H}^p_{D}(M).
$$
Therefore, as in the classical case, an abolute (weighted) cohomology class is uniquely represented by a $p$-form $\omega$ such that
$$\begin{cases}
    d\omega=\delta_f \omega =0 & \text{on } M\\
    \nu \lrcorner \, \omega =0 & \text{on } \partial M.
\end{cases}$$
A relative (weighted) cohomology class is uniquely represented by a $p$-form $\omega$ such that
$$\begin{cases}
    d\omega=\delta_f \omega =0 & \text{on } M\\
    \iota^* \omega =0 & \text{on } \partial M.
\end{cases}$$
\begin{remark}
    Note that the fact that the weighted cohomology group is isomorphic to the unweighted one  is well-known in weighted manifolds without boundary, see e.g. \cite[Id. (2.13)]{Lott03}.
\end{remark}
\begin{prop}
\label{vanishing_relative_cohomology}
    Let $(M,g, d\mu_f)$ be a compact Bakry--\'Emery manifold with smooth boundary $\partial M$. Assume that $\ric_f^{(p)} \ge 0$ for some $p \in \{1, \dots ,n\}$ and $\inf_{\partial M} (\sigma_{n-p} +f_\nu) >0$. Thus $\mathcal{H}^p_{D}(M)=0$.
 \end{prop}
 {\it Proof.} Let $\omega$ be a $p$-form satisfying $d\omega=\delta_f \omega=0$ on $M$ and $\iota^* \omega=0$ on $\partial M$. By the Reilly formula applied to $\omega$, we get
\begin{eqnarray*}
    0&=&\int_M |\nabla \omega|^2 \, d\mu_f + \int_M \langle \ric^{(p)}_f \omega , \omega \rangle \, d\mu_f +\int_{\partial M} \langle S^{[n-p]}(\iota^* \star \omega), \iota^* \star \omega \rangle \, d \mu_f + \int_{\partial M}f_\nu | \nu \lrcorner \, \omega |^2 \, d \mu_f\\
    &\ge &\int_M |\nabla \omega|^2 \, d\mu_f + \int_{\partial M} (\sigma_{n-p} +f_\nu)|\nu \lrcorner \,  \omega|^2 \, d\mu_f\\
    &\ge &\int_M |\nabla \omega|^2 \, d\mu_f + \inf_{\partial M} (\sigma_{n-p} +f_\nu) \int_{\partial M} |\nu \lrcorner \,  \omega|^2 \, d\mu_f. 
\end{eqnarray*}
Therefore, we deduce that the form $\omega$ is parallel and that $\nu \lrcorner \, \omega =0$, which implies that $\omega$ vanishes on $M$. \hfill$\square$
 \begin{remark}
 Under the condition that the $p$-Ricci tensor $\ric_f^{(p)} \ge 0$ for some $p \in \{1, \dots ,n\}$, we can see by using the Reilly formula that if $\sigma_p(\partial M) \ge 0$ (resp. $\inf_{\partial M} (\sigma_{n-p} +f_\nu) \ge 0$) then any form $\omega$ such that $d\omega=\delta_f \omega =0$ and $\nu \lrcorner \, \omega=0$ on $\partial M$ (resp. $\iota^* \omega =0$ on $\partial M$) is necessarily parallel and satisfies $\nabla f \lrcorner \, \omega=0$. Therefore, $\delta \omega$ also vanishes.  
 \end{remark}   
We now derive a Poincar\'e-type inequality for differential forms in a Bakry--\'Emery manifold. 
\begin{thm}[Weighted Poincar\'e inequality]
\label{Poincare}
Let $(M^n,g,d\mu_f)$ be a compact weighted manifold with smooth boundary $\partial M$ and let $\nu$ be the inward unit normal vector field along the boundary.
    Assume that there exists $p \in \{1, \dots ,n-1\}$, $p \ne n/2$ such that $\ric_{f}^{(p)} \ge 0$ and $\sigma_{n-p}+f_\nu >0$ along the boundary. Then for any exact $p$-form $\omega$ on $\partial M$, we have
    \begin{equation}
    \label{poincare_inequality}
    \int_{\partial M} \langle S^{[p]}\omega, \omega \rangle  \, d \mu_f \le \int_{\partial M}\dfrac{1}{\sigma_{n-p}+f_\nu}|\delta_f^{\partial M} \omega|^2\, d \mu_f.\end{equation}
    Equality is achieved in \eqref{poincare_inequality} for a non-zero form $\omega$, if and only if the $p$-form $\hat{\omega}$ on M -- which satisfies $d \hat{\omega}=\delta_f \hat{\omega}=0$ on $M$ and $\iota^* \hat{\omega} = \omega$ on $\partial M$ -- is parallel, and $\delta_f^{\partial M} \omega=-(\sigma_{n-p}+f_\nu) \nu \lrcorner \, \hat{\omega}$. \\
    Furthermore, in the equality case, the following identities must hold along $\partial M$: $S^{[p]} (\omega) = -d^{\partial M}(\nu \lrcorner \,  \hat{\omega})$, $S^{[p-1]}(\nu \lrcorner \, \hat{\omega})=((n-1)H-\sigma_{n-p}) \, \nu \lrcorner \, \hat{\omega}$.
\end{thm}
{\it Proof.}
Let $\omega =d^{\partial M} \alpha$ be an exact $p$-form on the boundary with $\alpha \in \Omega^{p-1}(\partial M)$.  There exists a $(p-1)$-form $\hat{\alpha}$ such that $\delta_f d \hat{\alpha}=0$ and $\iota^* \hat{\alpha}=\alpha$. To see this, notice that, for the conformal metric $\overline{g}:=e^{2u}g$ with $u:=-\frac{f}{n-2p}$  on $M$, we have $\overline{\delta}=\delta^{\overline{g}}=e^{-2u}\cdot(\delta-(n-2p)\nabla u\lrcorner \,\,)=e^{\frac{2f}{n-2p}}\cdot(\delta+\nabla f\lrcorner \,\,)=e^{\frac{2f}{n-2p}}\cdot\delta_f$ on $p$-forms, so that $\delta_fd\hat{\alpha}=0$ is equivalent to $\overline{\delta}d\hat{\alpha}=0$ on $(M,\overline{g})$, therefore \cite[Thm. 2]{DuffSpencer52} applies and yields the existence of $\hat{\alpha}$.
 
We set $\hat{\omega} := d \hat{\alpha} \in \Omega^p(M)$. The form $\hat{\omega}$ satisfies
$$\begin{cases}
    d\hat{\omega}=0 & \text{on } M \\
    \delta_f\hat{\omega}=0 & \text{on } M \\
    \iota^* \hat{\omega} =\omega & \text{on } \partial M.
\end{cases}$$
Note that the form $\hat{\omega}$ is unique due to Proposition \ref{vanishing_relative_cohomology}. By the Reilly formula \eqref{WReilly_formula} applied to $\hat{\omega}$, and given that $\ric_f^{(p)} \ge 0$, we get
\begin{eqnarray}
\label{eq:est_term_Reilly}
    0 &\ge& 2\int_{\partial M} \langle \delta^{\partial M}_f \omega , \nu \lrcorner \, \hat{\omega} \rangle \, d \mu_f  + \int_{\partial M} \mathcal{B}_f(\hat{\omega}, \hat{\omega}) \, d \mu_f \nonumber\\
    &\ge & 2\int_{\partial M} \langle \delta^{\partial M}_f \omega , \nu \lrcorner \, \hat{\omega} \rangle \, d \mu_f  +  \int_{\partial M} (\langle S^{[p]}\omega, \omega \rangle + (\sigma_{n-p} +f_\nu) |\nu \lrcorner \, \hat{\omega}|^2) \, d \mu_f.
\end{eqnarray}
Now, by using the condition $\sigma_{n-p} + f_\nu >0$, we have
\begin{eqnarray}
\label{eq:est_term_Reilly2}
2 \langle \delta^{\partial M}_f \omega , \nu \lrcorner \, \hat{\omega} \rangle  &=&\left|  (\sigma_{n-p}+f_\nu)^{1/2} \nu \lrcorner \, \hat{\omega} + \dfrac{1}{(\sigma_{n-p}+f_n)^{1/2}} \delta_f^{\partial M}\omega  \right|^2 \nonumber\\
&&-(\sigma_{n-p}+f_\nu) |\nu\lrcorner \, \hat{\omega}|^2 -\dfrac{1}{\sigma_{n-p}+f_\nu}|\delta_f^{\partial M} \omega|^2 \nonumber\\
&\ge & -(\sigma_{n-p} +f_\nu) |\nu\lrcorner \, \hat{\omega}|^2 -\dfrac{1}{\sigma_{n-p}+f_\nu}|\delta_f^{\partial M} \omega|^2,
\end{eqnarray}
from which we obtain
\begin{equation*}
    2 \langle \delta^{\partial M}_f \omega , \nu \lrcorner \, \hat{\omega} \rangle + (\sigma_{n-p} +f_\nu) |\nu\lrcorner \, \hat{\omega}|^2 \ge -\dfrac{1}{\sigma_{n-p}+f_\nu}|\delta_f^{\partial M} \omega|^2.
\end{equation*}
Combining the previous inequality with \eqref{eq:est_term_Reilly}, we deduce
$$0 \ge \int_{\partial M} (\langle S^{[p]}\omega, \omega \rangle  -\dfrac{1}{\sigma_{n-p}+f_\nu}|\delta_f^{\partial M} \omega|^2) \, d \mu_f,$$
which is the stated inequality. Assume now that \eqref{poincare_inequality} is an equality. In that case, equality must hold in both \eqref{eq:est_term_Reilly} and \eqref{eq:est_term_Reilly2}, from which it follows that $\hat{\omega}$ is parallel, $S^{[n-p]}(\iota^* (\star \hat{\omega}))=\sigma_{n-p} \, \iota^* (\star \hat{\omega})$ (that is $S^{[p-1]}(\nu \lrcorner \, \hat{\omega})=((n-1)H-\sigma_{n-p}) \, \nu \lrcorner \, \hat{\omega}$), and $\delta_f^{\partial M} \omega=-(\sigma_{n-p}+f_\nu) \nu \lrcorner \, \hat{\omega}$. 
Conversely, if we assume $\hat{\omega}$ to be parallel, then,  by \eqref{ddeltabord}, we have $S^{[p]} (\omega) = -d^{\partial M}(\nu \lrcorner \,  \hat{\omega})$. Therefore, using $\delta_f^{\partial M} \omega=-(\sigma_{n-p}+f_\nu) \nu \lrcorner \, \hat{\omega}$, we obtain the following:
\begin{eqnarray*}
    \int_{\partial M} \< S^{[p]} \omega , \omega \> \, d \mu_f &=& -\int_{\partial M} \< d^{\partial M}(\nu \lrcorner \,  \hat{\omega}) , \omega \> \, d \mu_f \\
    &=& -\int_{\partial M} \< \nu \lrcorner \,  \hat{\omega} , \delta_f^{\partial M} \omega \> \, d \mu_f \\
    &=& \int_{\partial M}  \dfrac{|\delta_f^{\partial M} \omega|^2}{\sigma_{n-p}+f_\nu}  \, d \mu_f.
\end{eqnarray*}
This concludes the proof.  \hfill$\square$ 
\begin{remark} \label{eq_case_unweighted}
    We assume $\sigma_{n-p} >0$. If the equality case is attained in \eqref{poincare_inequality} for an exact $p$-form $\omega$, then we have the equality in the unweighted case \cite[Thm. 3.1]{GHR}, namely:
    \begin{equation}
    \label{un_poinc_eq_case}\int_{\partial M} \langle S^{[p]}\omega, \omega \rangle  \, d \mu_g = \int_{\partial M}\dfrac{1}{\sigma_{n-p}}|\delta^{\partial M} \omega|^2\, d \mu_g.    
    \end{equation}
    Indeed, if the form $\hat{\omega}$ is parallel on M and satisfies $\iota^* \hat{\omega} = \omega$, $S^{[p]} (\omega) = -d^{\partial M}(\nu \lrcorner \,  \hat{\omega})$, $S^{[p-1]}(\nu \lrcorner \, \hat{\omega})=((n-1)H-\sigma_{n-p}) \, \nu \lrcorner \, \hat{\omega}$, and $\delta_f^{\partial M} \omega=-(\sigma_{n-p}+f_\nu) \nu \lrcorner \, \hat{\omega}$ on $\partial M$. By \eqref{ddeltabord}, we also have  $\delta^{\partial M} \omega=-\sigma_{n-p} \, \nu \lrcorner \, \hat{\omega}$, so we immediately obtain \eqref{un_poinc_eq_case}.
\end{remark}
Our next result establishes a Poincar\'e-type inequality for differential forms in weighted manifolds in the case where the Bakry--\'Emery Ricci tensor is bounded below by a real number $c \ne 0$. 
\begin{thm}
\label{poincare_gen}
Let $(M^n,g,d\mu_f)$ be a compact weighted manifold with smooth boundary $\partial M$ and let $\nu$ be the inward unit normal vector field along the boundary.
    Assume that there exists $p \in \{1, \dots ,n-1\}$, $p \ne n/2$ such that $\ric_{f}^{(p)} \ge c \, p(n-p)$ for some $c  \ne 0$ and $\sigma_{n-p}+f_\nu >0$ along the boundary. Then for any exact $p$-form $\omega=d^{\partial M}\alpha$ on the boundary
    \begin{equation} \label{poin_gen_ricfpos}\int_{\partial M} \langle S^{[p]}\omega, \omega \rangle  \, d \mu_f \le \int_{\partial M}\dfrac{|\delta_f^{\partial M} \omega  -\frac{c  \, p(n-p)}{2} \alpha|^2}{\sigma_{n-p}+f_\nu}\, d \mu_f.\end{equation}
    Moreover, equality is achieved in \eqref{poin_gen_ricfpos} if and only if $\alpha=0$.
\end{thm}
{\it Proof.}
Let $\omega =d^{\partial M} \alpha$ be an exact $p$-form on the boundary with $\alpha \in \Omega^{p-1}(\partial M)$. We proceed as in Theorem \ref{Poincare}. There exists a $(p-1)$-form $\hat{\alpha}$ such that $\delta_f d \hat{\alpha}=0$ and $\iota^* \hat{\alpha}=\alpha$. We put $\hat{\omega} := d \hat{\alpha} \in \Omega^p(M)$. The form $\hat{\omega}$ satisfies
$$\begin{cases}
    d\hat{\omega}=0 & \text{on } M \\
    \delta_f\hat{\omega}=0 & \text{on } M \\
    \iota^* \hat{\omega} =\omega & \text{on } \partial M.
\end{cases}$$
By the Reilly formula applied to $\hat{\omega}$, and given that $\ric_f^{(p)} \ge c  \, p(n-p)$, and using the fact that $\int_M |\hat{\omega}|^2 \, d\mu_f =-\int_{\partial M} \langle \alpha,\nu \lrcorner \, \hat{\omega} \rangle \, d\mu_f$ by partial integration, we obtain
\begin{eqnarray}
\label{eq:est_terms_Reilly3}
    0 &\ge& -c  \, p(n-p)\int_{\partial M} \langle \alpha,\nu \lrcorner \, \hat{\omega} \rangle \, d\mu_f+2\int_{\partial M} \langle \delta^{\partial M}_f (\omega) , \nu \lrcorner \, \hat{\omega} \rangle \, d \mu_f  + \int_{\partial M} \mathcal{B}_f(\hat{\omega}, \hat{\omega}) d \mu_f \nonumber\\
    &\ge & 2\int_{\partial M} \langle \delta^{\partial M}_f \omega -\frac{c  \, p(n-p)}{2} \alpha, \nu \lrcorner \, \hat{\omega} \rangle \, d \mu_f  +  \int_{\partial M} (\langle S^{[p]}\omega, \omega \rangle + (\sigma_{n-p} +f_\nu) |\nu \lrcorner \, \hat{\omega}|^2) \, d \mu_f.
\end{eqnarray}
Now, by making use of the condition $\sigma_{n-p} + f_\nu >0$, we have
\begin{eqnarray*}2 \langle \delta^{\partial M}_f \omega  -\frac{c  \, p(n-p)}{2} \alpha, \nu \lrcorner \, \hat{\omega} \rangle  &=&\left|  (\sigma_{n-p}+f_\nu)^{1/2} \nu \lrcorner \, \hat{\omega} + \dfrac{1}{(\sigma_{n-p}+f_n)^{1/2}} \left(\delta_f^{\partial M}\omega  -\frac{c  \, p(n-p)}{2} \alpha \right)  \right|^2 \\
&&-(\sigma_{n-p}+f_\nu) |\nu\lrcorner \, \hat{\omega}|^2 -\dfrac{1}{\sigma_{n-p}+f_\nu}|\delta_f^{\partial M} \omega  -\frac{c  \, p(n-p)}{2} \alpha|^2 \\
&\ge & -(\sigma_{n-p} +f_\nu) |\nu\lrcorner \, \hat{\omega}|^2 -\dfrac{1}{\sigma_{n-p}+f_\nu}|\delta_f^{\partial M} \omega  -\frac{c  \, p(n-p)}{2} \alpha|^2.
\end{eqnarray*}
Therefore,
\begin{equation*}
    2 \langle \delta^{\partial M}_f \omega  -\frac{c  \, p(n-p)}{2} \alpha, \nu \lrcorner \, \hat{\omega} \rangle + (\sigma_{n-p} +f_\nu) |\nu\lrcorner \, \hat{\omega}|^2 \ge -\dfrac{1}{\sigma_{n-p}+f_\nu}|\delta_f^{\partial M} \omega  -\frac{c  \, p(n-p)}{2} \alpha|^2.
\end{equation*}
Combining this last inequality with \eqref{eq:est_terms_Reilly3} gives
$$0 \ge \int_{\partial M} (\langle S^{[p]}\omega, \omega \rangle -\dfrac{1}{\sigma_{n-p}+f_\nu}|\delta_f^{\partial M} \omega  -\frac{c  \, p(n-p)}{2} \alpha|^2) \, d \mu_f.$$
Assume now that \eqref{poin_gen_ricfpos} is an equality, then the form $\hat{\omega}$ is parallel,  
$S^{[p-1]}(\nu \lrcorner \, \hat{\omega})=((n-1)H-\sigma_{n-p}) \, \nu \lrcorner \, \hat{\omega}$, and 
$\delta_f^{\partial M} \omega=\frac{c  \, p(n-p)}{2} \alpha-(\sigma_{n-p}+f_\nu) \nu \lrcorner \, \hat{\omega}$.
Thus, combining Equation \eqref{deltabord_Weighted} with the previous identities yields the following:
\begin{eqnarray*}
  \delta_f^{\partial M} \omega&=& \iota^*(\delta_f \hat{\omega})+\nu \lrcorner \, \nabla_\nu \hat{\omega} +S^{[p-1]}(\nu \lrcorner \, \hat{\omega}) -(n-1)H \nu \lrcorner \, \hat{\omega} - f_\nu \, \nu \lrcorner \, \hat{\omega}\\
  &=& S^{[p-1]}(\nu \lrcorner \, \hat{\omega}) -(n-1)H \nu \lrcorner \, \hat{\omega} - f_\nu \, \nu \lrcorner \, \hat{\omega} \\
  &=& -(\sigma_{n-p} +f_\nu) \, \nu \lrcorner \, \hat{\omega} .
 \end{eqnarray*}
A comparison between this expression of $\delta_f^{\partial M} \omega$ and the identity  $\delta_f^{\partial M} \omega=\dfrac{c  \, p(n-p)}{2} \alpha-(\sigma_{n-p}+f_\nu) \nu \lrcorner \, \hat{\omega}$ gives $\alpha=0$.
 \hfill$\square$

In the next theorem, we shall consider the particular case where the boundary $\partial M$ of the manifold $M$ is isometrically immersed in $\R^{d}$, where $d \ge n-1$. Our result is an extension of \cite[Thm. 3.3]{GHR} to the case of weighted manifolds. The scalar case is studied in \cite[Thm. 1.2, (1)]{TuHuang2019}.
\begin{thm}
\label{boundary_isom_imm_eucliean}
Let $(M^n,g,d\mu_f)$ be a compact weighted manifold such that the boundary $\partial M$ of $M$ is isometrically immersed in $\R^{d}$.
    Assume that there exists $p \in \{1, \dots ,n-1\}$, $p \ne n/2$ such that $\ric_{f}^{(p)} \ge 0$ and $\sigma_{n-p}+f_\nu >0$ along the boundary. Then
    \begin{eqnarray}
    \label{ineq_imm_euclidean}
    \int_{\partial M} H \, d\mu_f \le \int_{\partial M} \dfrac{1}{\sigma_{n-p} +f_\nu} \left((n-1) |H_0|^2 - \dfrac{p-1}{(n-1)(n-2)} \scal^{\partial M}  + \dfrac{1}{n-1} |d^{\partial M} f|^2\right)\, d\mu_f. 
\end{eqnarray}
Here $H_0$ denotes the mean curvature of the immersion $\iota_0$ and $\scal^{\partial M}$ the scalar curvature of $\partial M$.
\\
When the function 
$f$ is constant and $M$ is the 
$n$-dimensional flat disk $\mathbb{D}^n$, standardly embedded in an 
$n$-dimensional affine subspace of 
$\R^d$, equality is achieved in \eqref{ineq_imm_euclidean}. Conversely, assume that \eqref{ineq_imm_euclidean} is an equality, $\partial M$ is connected, $p \ge 2$, and $\iota_0(\partial M)$ is contained in an $n$-dimensional affine subspace. Then, up to a rescaling of the metrics on $M$ and $\R^d$, the manifold $(M,g)$ must be isometric to the flat disk $\mathbb{D}^n$, standardly embedded in that subspace.
\end{thm}
{\it Proof.}
Let $x_1, \dots , x_{d}$ be the standard coordinates on $\R^{d}$. For $i=1, \dots, d$, the unit parallel vector field $\partial_{x_i}$ splits into $\partial_{x_i}=\partial_{x_i}^T +\partial_{x_i}^\perp\in T \partial M \oplus T^\perp \partial M$. For $i_1, \dots , i_p \in \{1, \dots , d\}$, let us denote by $I$ the multi-index $(i_1, \dots, i_p)$, $dx_I:= dx_{i_1} \wedge \dots \wedge dx_{i_p}$ and $\omega_I:= dx_{i_1}^T \wedge \dots \wedge dx_{i_p}^T=dx_I^T$, where $dx_{j}^T =\iota^*_0 dx_{j}$ and $dx_{I}^T =\iota^*_0 dx_{I}$ with $\iota_0 : \partial M \rightarrow \R^d$ being the immersion.  We denote by $(e_1, \dots , e_{n-1}, \nu_1, \dots , \nu_m)$ a pointwise o.n.b. of $\R^d$ adapted to the decomposition $T \partial M \oplus T^\perp \partial M$, where $m=d-n+1$. Recall that for any $p$-forms $\alpha$ and $\beta$ on $\partial M$,
$$\sum_{k=1}^{n-1} \< e_k \lrcorner \, \alpha, e_k \lrcorner\, \beta\>=p\< \alpha, \beta\>, $$
due to the formula $\sum_{k=1}^{n-1} e_k^* \wedge ( e_k \lrcorner \, \alpha)=p \alpha$, where $e_k^*$ is the 1-form dual to $e_k$. As a consequence, we immediately obtain:
\begin{equation}
\label{formula_inner_p}\sum_{i=1}^d \<d{x_i}^T \wedge \alpha ,d {x_i}^T \wedge \beta \> =\sum_{k=1}^{n-1} \< e_k^* \wedge \alpha, e_k^* \wedge \beta\>=(n-1-p)\< \alpha, \beta\>.     
\end{equation}
Now, by the Poincar\'e inequality \eqref{poincare_inequality} applied to the exact form $\omega_I$, we obtain:
\begin{equation}
\label{sum_ptuple}
\int_{\partial M} \langle S^{[p]}\omega_I, \omega_I \rangle  \, d \mu_f \le \int_{\partial M}\dfrac{1}{\sigma_{n-p}+f_\nu}|\delta_f^{\partial M} \omega_I|^2\, d \mu_f.
\end{equation}
We will follow \cite[Thm. 3.3]{GHR} to compute the sum over all multi-indices $I$ of the terms appearing in the previous inequality.
By \cite[Eq. (3.6)]{GHR}, we have for any symmetric endomorphism $A$ of $T \partial M$, $\sum_I \< A^{[p]}\omega_I, \omega_I \> = p!\binom{n-1}{p} \frac{p}{n-1} \textrm{tr} A$. 
As a consequence, we deduce:
\begin{equation}
\label{weingartenomegaI}
    \sum_I \< S^{[p]}\omega_I, \omega_I \> = p!\binom{n-1}{p} p H.
\end{equation}
 Now, we will compute the r.h.s. of \eqref{sum_ptuple}. We will show that:
 \begin{eqnarray}
\label{sum_deltaf}
    \sum_{I} |\delta_f^{\partial M} \omega_I|^2 =p!\binom{n-1}{p} p  \left((n-1) |H_0|^2 - \dfrac{p-1}{(n-1)(n-2)} \scal^{\partial M}  + \dfrac{1}{n-1} |d^{\partial M} f|^2\right).
\end{eqnarray}
Note that, once \eqref{sum_deltaf} is proven, the substitution of \eqref{weingartenomegaI} and \eqref{sum_deltaf} into \eqref{sum_ptuple} immediately implies the desired inequality \eqref{ineq_imm_euclidean}. 
Let us now prove \eqref{sum_deltaf}. Observe first that 
\begin{eqnarray}
\label{deltaf}
    \sum_{I} |\delta_f^{\partial M} \omega_I|^2 =\sum_I |\delta^{\partial M} \omega_I|^2 +\sum_I |\nabla^{\partial M} f \lrcorner \, \omega_I|^2 + 2\sum_I \langle \delta^{\partial M} \omega_I, \nabla^{\partial M} f \lrcorner \, \omega_I \rangle.
\end{eqnarray}
The first term is computed in \cite[Lem. 2.2]{Savo05}:
\begin{equation*}
    \sum_I |\delta^{\partial M} \omega_I|^2= p!\binom{n-1}{p} p \left((n-1) |H_0|^2 - \dfrac{p-1}{(n-1)(n-2)} \scal^{\partial M}  \right).
\end{equation*}
Next, we estimate the second term of Equation \eqref{deltaf}. 
\begin{eqnarray*}
    \label{eq:nablafimmersion}
     \sum_I |\nabla^{\partial M} f \lrcorner \, \omega_I|^2 &=& \sum_I \< \omega_I, d^{\partial M} f \wedge (\nabla^{\partial M} f \lrcorner \, \omega_I) \> \\
    &=&  |d^{\partial M} f|^2 \sum_I | \omega_I|^2 - \sum_I |d^{\partial M} f \wedge \omega_I|^2\\
    &\stackrel{\eqref{formula_inner_p}}{=}& |d^{\partial M} f|^2 ((n-1)(n-2) \dots (n-p)-(n-2)(n-3) \dots (n-p-1)) \\
    &=&  p! \binom{n-1}{p} \dfrac{p}{n-1} |d^{\partial M} f|^2.
\end{eqnarray*}

We still have to compute the term $2 \sum_I \< \delta^{\partial M} \omega_I , \nabla^{\partial M} f \lrcorner \omega_I \>$. We will show that this term vanishes. We denote by $\omega_{i_1\dots \widehat{i_k} \dots i_p}$ the $(p-1)$-form $dx_{i_1}^T \wedge \dots \widehat{dx_{i_k}^T}\wedge \dots \wedge dx_{i_p}^T$.
By \cite[Eq. (3.3)]{Savo05}, we have:
\begin{eqnarray*}
    \delta^{\partial M} \omega_I&=&\sum_{k=1}^{p}(-1)^{k+1}
    \II^{[p-1]}_{\partial_{x_{i_k}}^\perp}(\omega_{i_1 \dots \widehat{i_k} \dots i_p})
     -(n-1)\sum_{k=1}^{p}(-1)^{k+1} \< H_0, \partial_{x_{i_k}}^\perp \> \omega_{i_1 \dots \widehat{i_k} \dots i_p},
\end{eqnarray*}
and consequently
\begin{align*}
    \sum_I \< \delta^{\partial M} \omega_I , \nabla^{\partial M} f \lrcorner \, \omega_I \> &= 
    \sum_I \sum_{k=1}^{p} (-1)^{k+1} 
    \< d^{\partial M} f \wedge \II^{[p-1]}_{\partial_{x_{i_k}}^\perp}(\omega_{i_1 \dots \widehat{i_k} \dots i_p}), \omega_{i_1 \dots i_p} \> \\
    &-(n-1) \sum_I \sum_{k=1}^{p} (-1)^{k+1} \< H_0, \partial_{x_{i_k}}^\perp \>
    \< d^{\partial M} f\wedge\omega_{i_1 \dots \widehat{i_k} \dots i_p}, \omega_{i_1 \dots i_p} \>.
\end{align*}
It follows that:
\begin{eqnarray*}
   && \sum_I \< \delta^{\partial M} \omega_I , \nabla^{\partial M} f \lrcorner \, \omega_I \> \\
    &=& 
    (p-1)\sum_{j_1, \dots, j_{p-2}=1}^{d} \sum_{i,j=1}^{d}  
    \< d^{\partial M} f \wedge \II^{[p-1]}_{\partial_{x_{i}}^\perp}(dx_j^T)\wedge \omega_{j_1 \dots  j_{p-2}}), dx_i^T \wedge dx_j^T\wedge\omega_{j_1 \dots j_{p-2}} \> \\
    &&-p (n-1)\sum_{j_1, \dots , j_{p-1}} \sum_{i=1}^{d}  \< H_0, \partial_{x_{i}}^\perp \>
    \< d^{\partial M} f\wedge\omega_{j_1 \dots j_{p-1}}, dx_i^T\wedge\omega_{j_1 \dots j_{p-1}} \> \\
    &=& (p-1).(n-p)(n-p+1) \dots (n-3)  \sum_{i,j=1}^{d}  
    \< d^{\partial M} f \wedge \II_{\partial_{x_{i}}^\perp}(dx_j^T), dx_i^T \wedge dx_j^T \> \\
    &&-p (n-1). (n-p)(n-p+1) \dots (n-2) \sum_{i=1}^{d}  \< H_0, \partial_{x_{i}}^\perp \>
    \< d^{\partial M} f, dx_i^T \>.
\end{eqnarray*}
The last term is clearly zero. Let us now verify that 
$\sum_{i,j=1}^{d}  
    \< d^{\partial M} f \wedge \II_{\partial_{x_{i}}^\perp}(dx_j^T), dx_i^T \wedge dx_j^T \>$ also vanishes.  Indeed, by the fact that $\II_V^{[p]}(\alpha_1 \wedge \dots \wedge \alpha_p)= \sum_{j=1}^p \alpha_1 \wedge \dots \wedge \II_V(\alpha_j) \wedge \dots \wedge \alpha_p$ for each vector $V$ normal to $T \partial M$ and $1$-forms $\alpha_1, \dots \alpha_p$, we get:
\begin{eqnarray*}
    &&\sum_{i,j=1}^{d}  
    \< d^{\partial M} f \wedge \II_{\partial_{x_{i}}^\perp}(dx_j^T), dx_i^T \wedge dx_j^T \> \\&=& \sum_{i,j=1}^{d}\left(  
    \< d^{\partial M} f, dx_i^T\> \< \II_{\partial_{x_{i}}^\perp}(dx_j^T),  dx_j^T \> - \< d^{\partial M} f, dx_j^T\> \< \II_{\partial_{x_{i}}^\perp}(dx_j^T),  dx_i^T \> \right) \\
    &=&\sum_{i,j=1}^{d}\sum_{a=1}^m \<\partial_{x_{i}}^\perp, \nu_a\>
    \left(  
    \< d^{\partial M} f, dx_i\>  \<\II_{\nu_a}(dx_j^T),  dx_j^T \> - \< d^{\partial M} f, dx_j^T\> \< \II_{\nu_a}(dx_j^T),  dx_i^T \> 
    \right) \\
    &=&0.
\end{eqnarray*}
We have therefore shown that $\sum_I \< \delta^{\partial M} \omega_I , \nabla^{\partial M} f \lrcorner \, \omega_I \>=0$.

If the function 
$f$ is constant and $M$ is the 
$n$-dimensional flat disk $\mathbb{D}^n$, standardly embedded in some 
$n$-dimensional affine subspace of 
$\R^d$, it follows by direct computation that equality holds in \eqref{ineq_imm_euclidean}.
Conversely, assume that \eqref{ineq_imm_euclidean} is an equality, then, we must have for each $I$ that $\hat{\omega}_I$ is parallel, $S^{[p]} \omega_I= -d^{\partial M}(\nu \lrcorner \,  \hat{\omega}_I)$, $S^{[p-1]}(\nu \lrcorner \, \hat{\omega}_I)=((n-1)H-\sigma_{n-p}) \nu \lrcorner \, \hat{\omega}_I$, and $\delta_f^{\partial M} \omega_I=-(\sigma_{n-p}+f_\nu) \nu \lrcorner \, \hat{\omega}_I$. Additionally, $\nabla f \lrcorner \, \hat{\omega}_I=0$ for any $p$-tuple $I$. Furthermore, since $\sigma_{n-p}>0$, equality \eqref{un_poinc_eq_case} is guaranteed by Remark \ref{eq_case_unweighted}. 
Hence, under the conditions on $M$, the equality case in  \cite[Thm. 3.3]{GHR} shows that $M$ is the flat disk $\mathbb{D}^n$ and the parallel $p$-forms $\hat{\omega}_I$, indexed by all $p$-tuples
$I=(i_1, \dots , i_p)$ with $1 \le i_1 < \dots < i_p \le n$, form a basis for the space of $p$-forms on 
$M$. Since $\nabla f \lrcorner \, \hat{\omega}_I=0$ for every multi-index $I$, it follows that $\nabla f$ must be zero, and therefore $f$ is constant on $M$. 
\hfill$\square$
\\\\
The next theorem considers the particular case of a Bakry--\'Emery manifold whose boundary is isometrically immersed in the sphere $\mathbb{S}^d(c)$ of curvature $c$.
\begin{thm}
\label{immersion_boundary_sphere}
Let $(M^n,g,d\mu_f)$ be a compact weighted manifold with smooth boundary $\partial M$.
    Assume that there exists $p \in \{1, \dots ,n-1\}$, $p \ne n/2$ such that $\ric_{f}^{(p)} \ge c p(n-p)$  for some $c>0$ and $\sigma_{n-p}+f_\nu >0$ along the boundary. We also assume that $\partial M$ is isometrically immersed in $\mathbb{S}^{d}(c)$, and let $H_0$ be the mean curvature of the immersion $\iota_0$ and $\scal^{\partial M}$ the scalar curvature of $\partial M$. Then:
    \begin{eqnarray}
    \label{ineq_imm_sphere}
    \int_{\partial M} H \, d\mu_f < \int_{\partial M} \dfrac{(n-1) |H_0|^2 - \dfrac{p-1}{(n-1)(n-2)} \scal^{\partial M}  + \dfrac{1}{n-1} |d^{\partial M} f|^2 +c \, \frac{p-1 +p(n-p)}{4}}{\sigma_{n-p} +f_\nu} \, d\mu_f.  
\end{eqnarray}
\end{thm}
{\it Proof.} We proceed as in the proofs of Theorem \ref{boundary_isom_imm_eucliean} and \cite[Thm. 3.6]{GHR}. Let $\iota_1 : \partial M \rightarrow \mathbb{R}^{d+1}$ be the composition of the canonical embedding $\mathbb{S}^d(c) \rightarrow \R^{d+1}$ with the immersion $\iota_0 : \partial M \rightarrow \mathbb{S}^d(c)$. By inequality \eqref{poin_gen_ricfpos} applied to $\omega_I =d^{\partial M} \alpha_I$, where  $\alpha_I:=x_{i_1}|_{\partial M} \, dx_{i_2}^T \wedge \dots \wedge dx_{i_p}^T=\iota_1^* (x_{i_1} dx_{i_2} \wedge \dots \wedge dx_{i_p})$ and $I=(i_1, \dots ,i_p)$ is a $p$-tuple, we obtain
\begin{equation} \label{ineq_imm_sphere}\int_{\partial M} \langle S^{[p]}\omega_I, \omega_I \rangle  \, d \mu_f \le \int_{\partial M}\dfrac{|\delta_f^{\partial M} \omega_I  -\frac{c p(n-p)}{2} \alpha_I|^2}{\sigma_{n-p}+f_\nu}\, d \mu_f.\end{equation}
Now, we will sum over $I$. Given that $\omega_I=dx_I^T$, by \eqref{weingartenomegaI}, the l.h.s is  $\sum_I \< S^{[p]}\omega_I, \omega_I \> = p!\binom{n-1}{p} p H$. The sum of the r.h.s. of \eqref{ineq_imm_sphere} is
\begin{eqnarray*}
  \sum_I |\delta_f^{\partial M} \omega_I  -\frac{c p(n-p)}{2} \alpha_I|^2 &=&\sum_I  \left( |\delta_f^{\partial M} \omega_I|^2  +\frac{c^2 \, p^2(n-p)^2}{4} |\alpha_I|^2  -c p(n-p) \<\delta_f^{\partial M} \omega_I , \alpha_I \> \right). 
\end{eqnarray*}
First, by \eqref{sum_deltaf}, and denoting the mean curvature of the immersion $\iota_1$ by 
$H_1$, we obtain:
\begin{eqnarray*}
    \sum_{I} |\delta_f^{\partial M} \omega_I|^2 &=&p!\binom{n-1}{p} p  \left((n-1) |H_1|^2 - \dfrac{p-1}{(n-1)(n-2)} \scal^{\partial M}  + \dfrac{1}{n-1} |d^{\partial M} f|^2\right)\\
    &=& p!\binom{n-1}{p} p  \left((n-1) (|H_0|^2+c^2) - \dfrac{p-1}{(n-1)(n-2)} \scal^{\partial M}  + \dfrac{1}{n-1} |d^{\partial M} f|^2\right).
\end{eqnarray*}
Here we used the identity $|H_1|^2=|H_0|^2 + c^2$ that relates the mean curvature $H_1$ of $\iota_1$ and the mean curvature $H_0$ of $\iota_0$.  
Furthermore, from \cite{GHR}, $\sum_I |\alpha_I|^2=\dfrac{1}{c} (p-1)!\binom{n-1}{p-1}$ and $\sum_I \<\delta^{\partial M} \omega_I , \alpha_I \> = (p-1)!\binom{n-1}{p-1} (n-p)$. To conclude the computation of the sum of the r.h.s. of \eqref{ineq_imm_sphere}, it remains to determine $\sum_I \<\nabla^{\partial M} f \lrcorner\, \omega_I , \alpha_I \>$. Observe that $\sum_I \<\nabla^{\partial M} f \lrcorner\, \omega_I , \alpha_I \> =\sum_I \< \omega_I , d^{\partial M} f \wedge\alpha_I \> $, therefore
\begin{eqnarray*}
    \sum_I \<\nabla^{\partial M} f \lrcorner\, \omega_I , \alpha_I \> &=& 
    \sum_{i_1, \dots i_p=1}^{d+1} \< x_{i_1} dx_{i_1}^T \wedge dx_{i_2}^T \wedge \dots \wedge dx_{i_p}^T, d^{\partial M} f \wedge dx_{i_2}^T \wedge \dots \wedge dx_{i_p}^T \>.
\end{eqnarray*}
Note that the vector field $\sum_{i=1}^{d+1}x_i \partial_{x_i}$ is normal to $\partial M$ because it is a normal vector field for the immersion $\mathbb{S}^d(c) \rightarrow \R^{d+1}$, and therefore we deduce that $\sum_{i=1}^{d+1} \< x_i  dx_i, d^{\partial M} f \>=0$, from which it follows that $\sum_I \<\nabla^{\partial M} f \lrcorner\, \omega_I , \alpha_I \> =0$.  
From the previous computations, we can then conclude that:
\begin{align*}
   \sum_I |\delta_f^{\partial M} \omega_I  -\frac{c p(n-p)}{2} \alpha_I|^2 =&  
       p!\binom{n-1}{p} p  \left((n-1) |H_0|^2 - \dfrac{p-1}{(n-1)(n-2)} \scal^{\partial M}  + \dfrac{1}{n-1} |d^{\partial M} f|^2\right) \\
   &+ p!\binom{n-1}{p} p  c \left(n-1 + \dfrac{p(n-p)}{4}-(n-p)\right)\\
   =& p!\binom{n-1}{p} p  \left((n-1) |H_0|^2 - \dfrac{p-1}{(n-1)(n-2)} \scal^{\partial M}  + \dfrac{1}{n-1} |d^{\partial M} f|^2\right) \\
   &+ p!\binom{n-1}{p} p  c \left(p-1 + \dfrac{p(n-p)}{4}\right).
\end{align*}
Therefore
\begin{eqnarray*}
    \int_{\partial M} H \, d\mu_f \le \int_{\partial M} \dfrac{(n-1) |H_0|^2 - \dfrac{p-1}{(n-1)(n-2)} \scal^{\partial M}  + \dfrac{1}{n-1} |d^{\partial M} f|^2 +c \, \frac{p-1 +p(n-p)}{4}}{\sigma_{n-p} +f_\nu} \, d\mu_f. 
\end{eqnarray*}
If equality is achieved in the last inequality, then $\alpha_I$ must be 0 by Theorem \ref{poincare_gen}. On the other hand, the forms $\alpha_I$ span $\wedge^{p-1} T^* \partial M$ pointwise, leading to a contradiction. Thus, the inequality is necessarily strict, which concludes the proof of Theorem  \ref{immersion_boundary_sphere}. \hfill$\square$ 
\begin{remark}
    When the function $f$ is constant, we recover \cite[Thm. 3.6]{GHR}. The case $p=1$ (with $c=1$) was studied by Tu and Wang \cite[Ineq. (1.7)]{TuHuang2019} under the assumption that the $m$-dimensional Bakry--\'Emery tensor is bounded from below by $(m-1)g$.
\end{remark}
Further applications of Poincar\'e-type inequalities include an estimate of the first eigenvalue of the weighted boundary Laplacian restricted to exact $p$-forms. \\\\
First, let us recall that, as in the classical case, the weighted boundary Laplacian $\Delta_f^{\partial M}$ commutes with $d^{\partial M}$ and $\delta_f^{\partial M}$, which implies that it preserves the spaces of exact and co-exact forms. Denote by $\lambda_{1,p,f}'(\partial M)$ (resp. $\lambda_{1,p,f}''(\partial M)$)
the first eigenvalue of $\Delta_f^{\partial M}$ restricted to exact (resp. co-exact) 
$p$-forms. By the Hodge decomposition theorem, the first positive eigenvalue of weighted boundary Laplacian is $\lambda_{1,p,f}(\partial M) =\min(\lambda_{1,p,f}'(\partial M), \lambda_{1,p,f}''(\partial M))$.
Moreover, by differentiating eigenforms, the first eigenvalue on co-exact 
$p$-forms satisfies
$$\lambda_{1,p,f}'(\partial M)=\lambda_{1,p-1,f}''(\partial M).$$

\begin{prop}
\label{lower_bound_coexact_boundary}
Let $(M^n,g,d\mu_f)$ be a compact weighted manifold with smooth boundary $\partial M$ and let $\nu$ be the inward unit normal vector field along the boundary.
    Assume that there exists $p \in \{1, \dots ,n-1\}$, $p \ne n/2$ such that $\ric_{f}^{(p)} \ge 0$, $\sigma_p(\partial M) >0$ and $\inf_{\partial M} (\sigma_{n-p} +f_\nu )>0$. Then the first eigenvalue $\lambda_{1,p,f}'(\partial M)$ of the weighted boundary Laplacian acting on exact $p$-forms satisfies
  $$\lambda_{1,p,f}' (\partial M)\ge \sigma_p(\partial M)\, \inf_{\partial M}\big(\sigma_{n-p} + f_\nu \big). $$ 
\end{prop}
{\it Proof.}
Let $\alpha$ be a co-exact eigenform of $\Delta_f^{\partial M}$ associated with the eigenvalue $\lambda=\lambda_{1,p-1,f}''(\partial M)$.

Applying the weighted Poincar\'e inequality \eqref{poincare_inequality} to $\omega:=d^{\partial M} \alpha$, we obtain 
$$ \int_{\partial M} \langle S^{[p]}d ^{\partial M}\alpha, d ^{\partial M}\alpha \rangle  \, d \mu_f \le \int_{\partial M}\dfrac{1}{\sigma_{n-p}+f_\nu}|\Delta_f^{\partial M} \alpha|^2\, d \mu_f.$$
Hence, we get the following inequality:
$$ \sigma_p(\partial M)\int_{\partial M}  |d ^{\partial M}\alpha|^2  \, d \mu_f \le \dfrac{1}{\inf_{\partial M}\big(\sigma_{n-p}+f_\nu\big)}\int_{\partial M}|\Delta_f^{\partial M} \alpha|^2\, d \mu_f,$$
from which we deduce that:
$$\left( \frac{\lambda^2}{\inf_{\partial M}\big(\sigma_{n-p} + f_\nu \big)} -\sigma_p(\partial M) \lambda \right) \int_{\partial M}|\alpha|^2 d \mu_f \ge 0 ,$$
leading to
$$\lambda \ge \sigma_p(\partial M). \inf_{\partial M}\big(\sigma_{n-p} + f_\nu\big). $$
\hfill$\square$

\begin{remark}
\label{lower_bound_coexact_boundary_positive_ricci}
Under the assumptions of the previous proposition, and assuming in addition that 
$\ric_{f}^{(p)} \ge c \, p(n-p)$ for some $c  >0$, it follows that:
  $$\begin{array}{l}
     \lambda_{1,p,f}' (\partial M) \notin\left[\frac{C_p +c \, p(n-p) -\sqrt{C_p^2 +2 c \,  p(n-p) C_p}}{2},\frac{C_p +c \, p(n-p) +\sqrt{C_p^2 +2 c  \, p(n-p) C_p}}{2}\right], 
  \end{array}$$
where $C_p:=\sigma_p(\partial M). \inf_{\partial M}\big(\sigma_{n-p} + f_\nu\big)$. Indeed, let $\alpha$ be a co-exact eigenform of $\Delta_f^{\partial M}$ associated with the eigenvalue $\lambda=\lambda_{1,p-1,f}''(\partial M)$.
The inequality \eqref{poin_gen_ricfpos} applied to $\omega:=d^{\partial M} \alpha$, yields 
$$\int_{\partial M} \langle S^{[p]}d ^{\partial M}\alpha, d ^{\partial M}\alpha \rangle  \, d \mu_f  <  \int_{\partial M}\dfrac{|\Delta_f^{\partial M} \alpha  -\frac{c  \, p(n-p)}{2} \alpha|^2}{\sigma_{n-p}+f_\nu}\, d \mu_f.$$
Proceeding as in the proof of Proposition \ref{lower_bound_coexact_boundary} yields the inequality:
$$\left( \lambda^2 -\left(\sigma_p(\partial M). \inf_{\partial M}\big(\sigma_{n-p} + f_\nu \big) +c  \, p(n-p)\right)\lambda +  \frac{c ^2 \, p^2(n-p)^2}{4} \right) \int_{\partial M}|\alpha|^2 d \mu_f > 0 .$$
The claim follows from a straightforward  computation.
\end{remark}
\section{Some boundary value problems}
In this section, we shall extend some classical eigenvalue problems for differential forms to the setting of Bakry--\'Emery manifolds and, under suitable curvature assumptions, derive estimates for the first  eigenvalue. 

Recall that a $p$-form $\omega$ is called to be $f$-harmonic if $\Delta_f \omega=0$. 
We begin by stating the following result, which is an immediate consequence of \eqref{f-harmonic_fields_DN}.
\begin{prop}
\label{kernel_WDirichlet}
    Let $(M^n,g,d\mu_f)$ be a compact weighted manifold with smooth boundary $\partial M$.  Then there is no nontrivial $f$-harmonic $p$-form that vanishes along the boundary.
    \end{prop}
    In what follows, we introduce some eigenvalues problems on differential forms in Bakry--\'Emery manifolds: the Diriclet, Neumann, buckling and clamped plate problems. Note that the ellipticity of these problems continues to hold, as the principal symbol of each operator coincides with that of its classical counterpart. In the scalar case, such problems have been studied in \cite{MaDu2010, IliasShouman2019, BezerraXia2020}, among other references.
\\\\    
We first consider the following Dirichlet eigenvalue problem on differential forms in a weighted manifold $(M^n,g,d\mu_f)$ with boundary:
\begin{equation}\label{WDirichlet}
\left\{
\begin{array}{lll}
\Delta_f \omega &=\lambda \omega &\textrm{ on } M\\
\omega &=0&\textrm{ on }\partial M.
\end{array}
\right.
\end{equation}
Similarly to the unweighted case, this problem is elliptic in the sense of Lopatinski\u{\i}-Shapiro (see \cite[Def. 1.6.1]{Sc}). Therefore, there exists a Hilbert basis of the space of $L^2$-integrable $p$-forms on $(M^n,g,d\mu_f)$, with respect to the mesaure $d \mu_f$, consisting of smooth eigenforms of \eqref{WDirichlet} associated to an unbounded and positive sequence of eigenvalues $(\lambda_{i,p,f})_{i\geq 1}$. In addition, we have the following variational characterization of the first eigenvalue:
\begin{equation}
\label{var_char_WDir}
    \lambda_{1,p,f}=\inf\left\{\frac{ \int_M |d\omega|^2 \, d\mu_f+\int_M |\delta_f\omega|^2 \, d\mu_f}{\int_M |\omega|^2 \, d\mu_f} \mid\omega\in\Omega^p(M)\setminus\{0\},\;\omega_{|_{\partial M}}=0\right\}.
\end{equation}
Similarly, we consider the weighted version of the Neumann eigenvalue problem on $p$-forms:
\begin{equation}\label{WNeumann}
\left\{
\begin{array}{lll}
\Delta_f \omega &= \mu \omega &\textrm{ on } M\\
\nu \lrcorner \,\omega &=0&\textrm{ on }\partial M\\
\nu \lrcorner \,d\omega &=0&\textrm{ on }\partial M.
\end{array}
\right.
\end{equation}
This problem is well-posed too, and its first positive eigenvalue $\mu_{1,p,f}$ is characterized by:
\begin{equation}\label{var_char_WNeumann}
\mu_{1,p,f}=\inf\left\{\frac{ \int_M |d\omega|^2 \, d\mu_f+\int_M |\delta_f\omega|^2 \, d\mu_f}{\int_M |\omega|^2 \, d\mu_f} \mid \omega\in\Omega^p(M)\setminus\{0\},\;\nu \lrcorner \,\omega=0, \omega \in \mathcal{H}^p_{f,N}(M)^\perp\right\}.
\end{equation}
If the absolute cohomology space of 
$M$ is nontrivial, then the smallest eigenvalue of problem \eqref{WNeumann} is $\mu_{0,p,f}=0$, and the corresponding eigenspace is $\mathcal{H}^p_{f,N}(M)$.

Following the study of the buckling and the clamped plate operators carried out in \cite{EGHMR}, one can extend their definitions and properties to the context of this paper. We define the buckling eigenvalue problem acting on differential forms in weighted manifolds as follows: 
\begin{equation}\label{buckling_weighted_forms}
\left\{
\begin{array}{lll}
\Delta_f^2 \omega&=\Lambda \Delta_f \omega &\textrm{ on } M\\ \omega
&=0&\textrm{ on }\partial M\\
\nabla_\nu\omega&=0&\textrm{ on }\partial M ,\end{array}\right.
\end{equation}
for some real constant $\Lambda$. Likewise, we generalize the clamped plate eigenvalue problem in the following way:
\begin{equation}\label{clamped_plate_weighted_forms}
\left\{
\begin{array}{lll}\Delta_f^2 \omega&=\Gamma \omega &\textrm{ on } M\\
\omega &=0&\textrm{ on }\partial M\\
\nabla_\nu\omega&=0&\textrm{ on }\partial M,
\end{array}
\right.
\end{equation}
 for some real constant $\Gamma$. Equations \eqref{ddeltabord} and \eqref{deltabord_Weighted} show that problems \eqref{buckling_weighted_forms} and \eqref{clamped_plate_weighted_forms} have equivalent boundary conditions, namely $\omega_{|_{\partial
M}}=0$, $\nu \lrcorner \, d \omega =0$ and $\iota^* \delta_f \omega =0$ on $\partial M$. 
For the respective problems \eqref{buckling_weighted_forms} and \eqref{clamped_plate_weighted_forms}, there exist monotonously nonincreasing positive real sequences of eigenvalues of finite multiplicities, denoted $(\Lambda_{i,p,f})_{i \ge 1}$ and $(\Gamma_{i,p,f})_{i \ge 1}$ 
  respectively. The corresponding eigenforms are smooth and constitute a Hilbert basis for the space of $L^2$-integrable 
$p$-forms on $(M,g,d\mu_f)$. 
Furthermore, $\Lambda_{1,p,f}$ is characterized by
\begin{equation}\label{eq:charb}
\Lambda_{1,p,f}=\inf\left\{\frac{\int_M |\Delta_f\omega|^2 \, d\mu_f}{
\int_M |d\omega|^2 \, d\mu_f+\int_M|\delta_f \omega |^2 \, d\mu_f} \mid \omega\in\Omega^p(M)\setminus\{0\},\,\omega_{|_{\partial
M}}=0\textrm{ and }\nabla_\nu\omega_{|\partial M}=0\right\}.
\end{equation}
Note also that $\Gamma_{1,p,f}$ is characterized by
\begin{equation}\label{eq:charcp}
\Gamma_{1,p,f}=\inf\left\{\frac{\int_{M}|\Delta_f\omega|^2 \, d \mu_f}{
\int_M| \omega|^2 \, d \mu_f} \mid \omega\in\Omega^p(M)\setminus\{0\}\mid\omega_{|_{\partial
M}}=0\textrm{ and }\nabla_\nu\omega_{|\partial M}=0\right\}.
\end{equation}
In all the variational characterizations of the preceding problems, the infimum is in fact a minimum, achieved by any eigenform associated with the corresponding eigenvalue of the problem.
\begin{remark}
	\enumerate
	\item Note that the operator $\Delta_f$ does not commute with the Hodge star operator, so the above eigenvalue problems are not invariant under the Hodge star operator in general.
	\item   Many other eigenvalue problems, such as the Robin problem (see \cite{EGH}) and the biharmonic Steklov problem (see \cite{EGHM, AbouAssali2025}), can be generalized in the same way, as well as the variational characterization of their respective first eigenvalues.
\end{remark}
   We now present some properties of the eigenvalues of the aforementioned problems that remain valid in the weighted case. Actually, due to the variational characterizations mentioned above, these results, originally proved in the classical (unweighted) setting in \cite{EGHMR}, admit a straightforward generalization to our setting of Bakry--\'Emery manifolds. Furthermore, in the scalar case, corresponding to operators acting on functions, our general framework recovers a number of important results from \cite{IliasShouman2019}. 
 \begin{thm}\label{BuckCPD}
Let $(M^n,g,d\mu_f)$ be a compact weighted manifold with smooth boundary $\partial M$. We then have: 
\begin{enumerate}
    \item $\Gamma_{1,p,f}<\Lambda_{1,p,f}^2$ for $0\leq p\leq n$;
    \item $\Lambda_{1,p,f} \lambda_{1,p,f}< \Gamma_{1,p,f}$ for $0\leq p\leq n$;
    \item $\lambda_{1,p,f}<\sqrt{\Gamma_{1,p,f}}<\Lambda_{1,p,f}$ for $0\leq p\leq n$;
    \item $\inf(\lambda_{1,p+1,f}, \lambda_{1,p-1,f}) \le \Lambda_{1,p,f}$ for $1\leq p\leq n-1$;
\item $\lambda_{1,1,f}\leq \Lambda_{1,0,f}$ and $\lambda_{1,n-1,f}\leq \Lambda_{1,n,f}$.
\item $\mu_{1,p,f} \le \Lambda_{1,p,f}$ for $0 \le p \le n$. Moreover, if the absolute cohomology space $\mathcal{H}_{N}^p(M)$ is trivial, this inequality is strict.
\end{enumerate}
\end{thm}
Our next goal is to establish an analogue of \cite[Prop. 2.9]{BrandingHabib22} for manifolds with boundary. For any real number $N \in ]-\infty,0[\cup ]n-p+1, \infty[$, let
$\ric_{N,f}^{(p)} :=\ric_{f}^{(p)}-\dfrac{1}{N-(n-p+1)} (df\wedge (\nabla f \lrcorner \,))$. 
\begin{prop}
\label{WDestimate}
Let $(M^n,g,d\mu_f)$ be a compact weighted manifold with smooth boundary $\partial M$ and let $\nu$ be the inward unit normal vector field along the boundary.
    Assume that $\ric_{N,f}^{(p)} \ge c  p (n-p)$ for some $c  >0$, some $p \in \{1, \dots ,n-1\}$ and some $N \in ]-\infty,0[\cup ]n-p+1, \infty[$. Then
    $$\lambda_{1,p,f} > \frac{c  \, p(n-p)}{\max(\frac{p}{p+1}, \frac{N-1}{N})}.$$
\end{prop}
{\it Proof.} Let $\omega$ be an eigenform of the weighted Dirichlet problem \eqref{WDirichlet} associated with the first eigenvalue $\lambda_{1,p,f}$. Considering that $\omega$ vanishes on the boundary, we obtain by the weighted Reilly formula \eqref{WReilly_formula} applied to $\omega$:
\begin{equation}
    \label{eq1_wei_Dir}
\int_M (|d \omega|^2 +|\delta_f \omega|^2)\, d \mu_f =\int_M (|\nabla \omega|^2 +\langle \ric_f^{(p)}(\omega), \omega \rangle ) \, d\mu_f.
\end{equation}
We now use the pointwise inequality $|\nabla \omega|^2 \ge \frac{1}{p+1}|d\omega|^2 +\frac{1}{n-p+1} |\delta \omega|^2$ (see \cite[Lemme 6.8]{GM}), which yields
\begin{eqnarray*}
    |\nabla \omega |^2 &\ge & \frac{1}{p+1}|d\omega|^2 +\dfrac{1}{n-p+1}|\delta_f \omega - \nabla f \lrcorner \, \omega|^2.
\end{eqnarray*}
Using the fact that $|\delta_f \omega - \nabla f \lrcorner \, \omega|^2 \ge |\delta_f \omega|^2 - |\nabla f \lrcorner \, \omega|^2$ and observing the inequality $\dfrac{(a-b)^2}{s} \ge \dfrac{a^2}{N} -\dfrac{b^2}{N-s}$, which holds for all $a,b,N,s \in \R$ satisfying $N(N-s) >0$, we find  (see \cite{BrandingHabib22})
$$\dfrac{1}{n-p+1}|\delta_f \omega - \nabla f \lrcorner \, \omega|^2 \ge \dfrac{1}{N}|\delta_f \omega|^2-\dfrac{1}{N-(n-p+1)}|\nabla f \lrcorner \, \omega|^2.$$ 
This allows us to obtain the following inequality, valid for any differential form $\omega$:
\begin{equation}
\label{estim_nabla_ricp}
   |\nabla \omega|^2 +\langle \ric_f^{(p)}(\omega), \omega \rangle ) \ge \dfrac{1}{p+1}|d \omega|^2 + \dfrac{1}{N}|\delta_f \omega|^2 + c  \, p(n-p) |\omega|^2. 
\end{equation}
Thus, substituting into \eqref{eq1_wei_Dir}, we obtain 
\begin{eqnarray*}
     \label{eq2_wei_Dir}
 \int_M ( |d \omega|^2 +|\delta_f \omega|^2)\, d \mu_f  &\ge& 
 \int_M \left(\dfrac{1}{p+1}|d \omega|^2 + \dfrac{1}{N}|\delta_f \omega|^2 + c  \, p(n-p) |\omega|^2\right)\, d \mu_f,
\end{eqnarray*}
which leads to
$$\dfrac{p}{p+1}\int_M |d\omega|^2 \, d\mu_f +\dfrac{N-1}{N}\int_M |\delta_f \omega|^2 \, d\mu_f \ge c  \, p(n-p)\int_M |\omega|^2 \, d \mu_f.$$
From this last inequality and the variational characterization \eqref{var_char_WDir}, it follows that
$$\lambda_{1,p,f} \ge \frac{c  \, p(n-p)}{\max(\frac{p}{p+1}, \frac{N-1}{N})}.$$
Furthermore, if equality holds, then $\omega$ is a conformal Killing $p$-form that vanishes along the boundary. Therefore, according to \cite[Lem. 3.17]{EGHMR},  $\omega$ satisfies $\nabla_\nu\omega_{|\partial M}=0$ along $\partial M$. We can then use the unique continuation theorem for elliptic second-order linear operators (see e.g. \cite[Thm. 1.4]{Salolectnotes2014} and \cite[Chap. VIII]{HoermanderlinearI}) to deduce that $\omega$ vanishes identically on $M$. This shows that the inequality must be strict.
\hfill$\square$ \\
\\
Combining statement 3 of Theorem \ref{BuckCPD} with Proposition \ref{WDestimate} yields the following:
\begin{cor}
    Under the same assumptions as in Proposition \ref{WDestimate}, one has:
    $$\Lambda_{1,p,f} > \frac{c  \, p(n-p)}{\max(\frac{p}{p+1}, \frac{N-1}{N})}$$
    and
    $$\Gamma_{1,p,f} > \left(\frac{c  \, p(n-p)}{\max(\frac{p}{p+1}, \frac{N-1}{N})}\right)^2.$$
\end{cor}
\begin{prop}
  Under the same conditions as in Proposition \ref{WDestimate}, and assuming additionally that $S^{[p]} \ge 0$, then 
\begin{equation}
\label{WN_eig_est}
    \mu_{1,p,f} \ge \frac{c  \, p(n-p)}{\max(\frac{p}{p+1}, \frac{N-1}{N})}.
\end{equation}
Moreover, if $\sigma_p >0$, then the above inequality is strict.   
\end{prop} 
{\it Proof.} We proceed as in the proof of Proposition \ref{WDestimate}. Let $\omega$ be an eigenform of the weighted Neumann problem \eqref{WNeumann} associated with the first eigenvalue $\mu_{1,p,f}$. By the weighted Reilly formula \eqref{WReilly_formula} applied to $\omega$, we obtain:
\begin{equation}
    \label{eq1_wei_Neumann}
\int_M (|d \omega|^2 +|\delta_f \omega|^2)\, d \mu_f =\int_M (|\nabla \omega|^2 +\langle \ric_f^{(p)}(\omega), \omega \rangle ) \, d\mu_f + \int_{\partial M} \langle S^{[p]} \iota^* \omega, \iota^* \omega \rangle \, d \mu_f. 
\end{equation}
Applying \eqref{estim_nabla_ricp} and the assumption $S^{[p]} \ge 0$ again yield
$$\max\left(\dfrac{p}{p+1}, \dfrac{N-1}{N} \right)\int_M (|d\omega|^2 + |\delta_f \omega|^2 )\, d\mu_f \ge c  \, p(n-p)\int_M |\omega|^2 \, d \mu_f.$$
The variational characterization \eqref{var_char_WNeumann} then allows us to deduce the announced estimate. \\
If $\sigma_p>0$ the estimate is clearly true because $S^{[p]} \ge 0$. In this case, we can further note that Equation \eqref{eq1_wei_Neumann} toghether with \eqref{estim_nabla_ricp} gives 
\begin{eqnarray*}
    \max\left(\dfrac{p}{p+1}, \dfrac{N-1}{N} \right)\int_M (|d\omega|^2 + |\delta_f \omega|^2 )\, d\mu_f &\ge& c  \, p(n-p)\int_M |\omega|^2 \, d \mu_f + \sigma_p(\partial M)\int_{\partial M} |\iota^* \omega|^2 \, d\mu_f \\
    &\ge& c  \, p(n-p)\int_M |\omega|^2 \, d \mu_f.
\end{eqnarray*}
Assume that equality holds in \eqref{WN_eig_est}. This implies that $\omega$ is a Killing form. Furthermore, since $\sigma_p >0$, we also have  $\iota^* \omega =0$, which means $\omega$ vanishes along the boundary. Consequently, by \cite[Lem. 3.17]{EGHMR}, $\omega$ must be zero. This completes the proof of the second part. \hfill$\square$

\end{document}